\documentclass[twoside]{article} 
\usepackage{graphicx}
\date{} 
\title{Asymptotics of some generalised sine-integrals}
\author{\sc R. B.\ Paris \\
{\em Division of Computing and Mathematics,} \\
{\em Abertay University, Dundee DD1 1HG, UK}}
\begin{document}
\def\f#1#2{\mbox{${\textstyle \frac{#1}{#2}}$}}
\def\dfrac#1#2{\displaystyle{\frac{#1}{#2}}}
\def\boldal{\mbox{\boldmath $\alpha$}}
\newcommand{\bee}{\begin{equation}}
\newcommand{\ee}{\end{equation}}
\newcommand{\lam}{\lambda}
\newcommand{\ka}{\kappa}
\newcommand{\al}{\alpha}
\newcommand{\la}{\lambda}
\newcommand{\ga}{\gamma}
\newcommand{\eps}{\epsilon}
\newcommand{\fr}{\frac{1}{2}}
\newcommand{\fs}{\f{1}{2}}
\newcommand{\g}{\Gamma}
\newcommand{\br}{\biggr}
\newcommand{\bl}{\biggl}
\newcommand{\ra}{\rightarrow}
\newcommand{\gtwid}{\raisebox{-.8ex}{\mbox{$\stackrel{\textstyle >}{\sim}$}}}
\newcommand{\ltwid}{\raisebox{-.8ex}{\mbox{$\stackrel{\textstyle <}{\sim}$}}}
\renewcommand{\topfraction}{0.9}
\renewcommand{\bottomfraction}{0.9}
\renewcommand{\textfraction}{0.05}
\newcommand{\mcol}{\multicolumn}
\date{}
\maketitle
\pagestyle{myheadings}
\markboth{\hfill \sc R. B.\ Paris  \hfill}
{\hfill \sc Asymptotic expansion of an integral\hfill}
\begin{abstract}
We obtain the asymptotic expansion for large integer $n$ of a generalised sine-integral
\[\int_0^\infty\left(\frac{\sin\,x}{x}\right)^{n}dx\]
by utilising the saddle-point method. This expansion is shown to agree with recent results of J. Schlage-Puchta in {\it Commun. Korean Math. Soc.} {\bf 35} (2020) 1193--1202 who used a different approach.

An asymptotic estimate is obtained for another related sine-integral also involving a large power $n$. Numerical results are given to illustrate the accuracy of this approximation. We also revisit the asymptotics of Ball's integral involving the Bessel function $J_\nu(x)$, which reduces to the above integral when $\nu=1/2$. 

\vspace{0.3cm}

\noindent {\bf Mathematics subject classification (2010):} 33E20, 34E05, 41A60
\vspace{0.1cm}
 
\noindent {\bf Keywords:} Sine integrals, asymptotic expansions, saddle-point approximation
\end{abstract}

\vspace{0.3cm}

\noindent $\,$\hrulefill $\,$

\vspace{0.3cm}
\begin{center}
{\bf 1.\ Introduction}
\end{center}
\setcounter{section}{1}
\setcounter{equation}{0}
\renewcommand{\theequation}{\arabic{section}.\arabic{equation}}
The expansion of the generalised sine-integral
\bee\label{e11}
I_n=\int_0^\infty\bl(\frac{\sin\,x}{x}\br)^{\!n}dx
\ee
for integer $n\to\infty$ has recently been considered by Schlage-Puchta \cite{SP}. However the method used seems to be unnecessarily involved and our aim here is to present a more direct computation using the well-known saddle-point method for Laplace-type integrals.
The interest in the integral $I_n$ stems from the fact that the intersection of the unit cube with a plane orthogonal to a diagonal and passing through the midpoint has $(n-1)$-measure equal to $2\sqrt{n}\,I_n/\pi$. These intersections arise naturally in certain probabilistic problems; see the references cited in \cite{SP}.

The second related sine-integral we consider is given by
\bee\label{e11a}
K_n=\int_0^\infty e^{-ax} \bl(1-\frac{\sin^2 x}{x^2}\br)^{\!n}dx\qquad (a>0)
\ee
for $n\to\infty$ when the parameter $a=O(1)$. An integral of this type was communicated to the author by H. Kaiser \cite{HK}. We employ a two-term saddle-point approximation to estimate the growth of $K_n$ for large $n$ and present numerical calculations to verify the accuracy of the resulting formula. In the final section we revisit the expansion for large $n$ of Ball's integral involving the Bessel function $J_\nu(x)$, which reduces to (\ref{e11}) when $\nu=\fs$.

\vspace{0.6cm}

\begin{center}
{\bf 2.\ The asymptotic expansion of $I_n$}
\end{center}
\setcounter{section}{2}
\setcounter{equation}{0}
\renewcommand{\theequation}{\arabic{section}.\arabic{equation}}
We begin by writing the integral in (\ref{e11}) as
\[I_n=\int_0^\pi\bl(\frac{\sin\,x}{x}\br)^{\!n}dx+R_n(x),\qquad R_n(x)=\int_\pi^\infty\bl(\frac{\sin\,x}{x}\br)^{\!n}dx.\]
It is easily seen that
\[|R_n(x)|<\int_\pi^\infty \frac{dx}{x^n}=\frac{\pi^{1-n}}{n-1}.\]
The remainder term $R_n(x)$ is therefore bounded by O($n^{-1}\pi^{-n}$)=O($n^{-1}e^{-n} (\pi/e)^{-n}$) and so is exponentially small as $n\to\infty$.

Let $\psi(x)=\log\,(x/\sin x)$, where $\psi(0)=0$ and $\psi(\pi)=\infty$. Then the integral over $[0,\pi]$ becomes
\[{\hat I}_n=\int_0^\pi\bl(\frac{\sin\,x}{x}\br)^{\!n}dx=\int_0^\pi e^{-n\psi(x)} dx.\]
This integral has a saddle point at $x=0$ and the integration path $[0,\pi]$ is the path of steepest descent through the saddle. If we now make the standard change of variable $\psi(x)=\tau^2$ discussed, for example, in \cite[p.~66]{ETC} we obtain
\[{\hat I}_n=\int_0^\infty e^{-n\tau^2}\, \frac{dx}{d\tau}\,d\tau.\]
From the expansion 
\[\tau^2=\log \bl(\frac{x}{\sin x}\br)=\frac{1}{6}x^2+\frac{1}{180}x^4+\frac{1}{2835}x^6+\frac{1}{37800}x^8+\frac{1}{467775} x^{10}+\cdots\  \]
valid for $|x|<\pi$,
we find by inversion of this series using {\it Mathematica} that
\[x=\sqrt{6}\bl\{\tau-\frac{1}{10}\tau^3
-\frac{13}{4200}\tau^5 + \frac{9}{14000}\tau^7 + \frac{17597}{77616000}\tau^9 + \frac
{4873}{218400000}\tau^{11}+\cdots\ \br\}\]
whence
\bee\label{e200}
\frac{dx}{d\tau}=\sqrt{6} \sum_{k=0}^\infty b_k \tau^{2k}\qquad (|\tau|<\tau_0).
\ee
The first few coefficients $b_k$ are
\[b_0=1,\quad b_1=-\frac{3}{10},\quad b_2=-\frac{13}{840},\quad b_3=\frac{9}{2000},\quad b_4=\frac{17597}{862400},
\]
\[b_5=\frac{53603}{218400000},\quad b_6= -\frac{124996631}{1629936000000},\quad b_7= -\frac{159706933}{4366252800000}, \ldots\ .\]
The circle of convergence of the series (\ref{e200}) is determined by the nearest point in the mapping $x\mapsto\tau$ where $dx/d\tau$ is singular (the point $x=\pi$ maps to $\infty$ in the $\tau$-plane); that is, where $\psi'(x)=0$. This yields $x=\tan\,x$, with the root $x_0\doteq1.43030\pi$, and hence
$\tau_0=|\log\,x_0/|\sin x_0|+\pi i|^{1/2}\doteq1.86894$.
Then we have
\[{\hat I}_n\sim \sqrt{6}\int_0^\infty e^{-n\tau^2} \sum_{k=0}^\infty b_k \tau^{2k}\,d\tau=\sqrt{\frac{3}{2n}}\sum_{k=0}^\infty \frac{b_k}{n^k} \int_0^\infty e^{-w} w^{k-1/2} dw\]
\bee\label{e201}
= \sqrt{\frac{3\pi}{2n}}\sum_{k=0}^\infty \frac{c_k}{n^k}\qquad(n\to\infty),
\ee
where the coefficients $c_k$ are defined by
\[c_k:=b_k\frac{\g(k+\fs)}{\g(\fs)}.\]
It follows that since we have extended the integration path in (\ref{e201}) beyond the circle of convergence of (\ref{e200}) the resulting asymptotic series is divergent.

\begin{table}[h]
\caption{\footnotesize{The coefficients $c_k$ for $1\leq k\leq12$ (with $c_0=1$).}}
\begin{center}
\begin{tabular}{|r|l|r|l|}
\hline
\mcol{1}{|c|}{$k$} & \mcol{1}{c|}{$c_k$} & \mcol{1}{c|}{$k$} & \mcol{1}{c|}{$c_k$} \\
[.1cm]\hline
&&&\\[-0.25cm]
1 & $-\f{3}{20}$ & 2 & $-\f{13}{1120}$ \\
&&&\\[-0.25cm]
3 & $+\f{27}{3200}$ & 4 & $+\f{527\,91}{394\,240\,0}$\\
&&&\\[-0.25cm]
5 & $+\f{482\,427}{665\,600\,00}$ & 6 & $-\f{124\,996\,631}{100\,352\,000\,00}$\\
&&&\\[-0.25cm]
7 & $-\f{527\,032\,878\,9}{13\,647\,872\,000\,0}$ & 8 & $-\f{747\,906\,350\,616\,1}{268\,461\,670\,400\,000}$\\
&&&\\[-0.25cm]
9 & $+\f{692\,197\,762\,461\,3}{565\,182\,464\,000\,00}$ & 10 & $+\f{107\,035\,304\,201\,928\,877\,41}{236\,585\,379\,430\,400\,000\,00}$\\
&&&\\[-0.25cm]
11& $+\f{509\,710\,579\,537\,397\,418\,9}{205\,726\,416\,896\,000\,000\,00}$ & 12 & $-\f{123\,979\,742\,078\,372\,360\,595\,39}
  {362\,078\,493\,736\,960\,000\,000\,0}$\\
[.1cm]\hline
\end{tabular}
\end{center}
\end{table}

Thus, neglecting exponentially small terms, we have the asymptotic expansion
\bee\label{e21}
I_n\sim \sqrt{\frac{3\pi}{2n}} \sum_{k=0}^\infty \frac{c_k}{n^k} \qquad (n\to\infty),
\ee
where the coefficients $c_k$ are listed in Table 1 for $0\leq k\leq 12$.
This expansion agrees with that obtained in \cite{SP} by less direct means, except for the value of the coefficient $c_{10}$.

An integral of a similar nature is
\[J_n=\int_0^\infty \bl(\frac{1-\cos x}{\fs x^2}\br)^{\!n}dx=\int_0^\infty \bl(\frac{\sin \fs x}{\fs x}\br)^{\!2n}dx=2I_{2n}.\]
From (\ref{e21}) its asymptotic expansion is therefore (to within exponentially small terms)
\[J_n\sim\sqrt{\frac{3\pi}{n}} \sum_{k=0}^\infty \frac{c_k}{(2n)^k}\qquad (n\to\infty).\]

\vspace{0.6cm}

\begin{center}
{\bf 3.\ An asymptotic estimate of another sine-integral}
\end{center}
\setcounter{section}{3}
\setcounter{equation}{0}
\renewcommand{\theequation}{\arabic{section}.\arabic{equation}}
In this section we consider the following integral
\bee\label{e31}
K_n=\int_0^\infty e^{-ax} \bl(1-\frac{\sin^2 x}{x^2}\br)^{\!n}dx\qquad (a>0)
\ee
for $n\to\infty$ (not necessarily an integer) when the parameter $a=O(1)$. 
We express $K_n$ as a Laplace-type integral in the form
\[K_n=\int_0^\infty e^{-n\psi(x)} f(x)\,dx,\]
where
\[\psi(x)=-\log \bl(1-\frac{\sin^2x}{x^2}\br),\qquad f(x)=e^{-ax}.\]
For large $n$ the exponential factor in the integrand consists of a series of peaks situated at $x=k\pi$, ($k=1, 2, \ldots$) of decreasing height controlled by the decay of $f(x)$; see Fig.~1 for a typical example. This is in marked contrast to the situation pertaining to the integral $I_n$ in (\ref{e11}), where the second and successive peaks are of height O($(k\pi)^{-n})$
($k\geq1$) and so are exponentially smaller than the (half) peak in $[0,\pi]$. 
\begin{figure}[th]
	\begin{center} \includegraphics[width=0.6\textwidth]{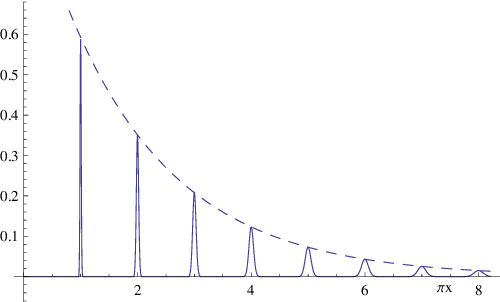}
\caption{\small{Plot of the integrand in (\ref{e31}) when $n=5000$ and $a=1/6$ with horizontal scale $\pi x$. The dashed curve represents $e^{-\pi ax}$.}} 
	\end{center}
\end{figure}
Routine calculations show that
\[\psi''(k\pi)=\frac{2}{(k\pi)^2},\quad \psi'''(k\pi)=-\frac{12}{(k\pi)^3},\quad \psi^{iv}(k\pi)=\frac{82}{(k\pi)^4}-\frac{8}{(k\pi)^2}.\]
Application of the two-term saddle-point approximation to the $k$th peak then yields the approximate contribution \cite[p.~48]{DLMF}, \cite[\S 1.2.3]{P}
\[2\sqrt{\frac{\pi}{2n\psi''(k\pi)}} \bl\{1+\frac{c_2}{n}\br\}e^{-k\pi a}
=k\pi\sqrt{\frac{\pi}{n}}\bl\{1+\frac{c_2}{n}\br\}e^{-k\pi a},\]
where
\[c_2=\frac{1}{2\psi''}\bl\{\frac{2f''}{f}-2\frac{\psi'''}{\psi''}\,\frac{f'}{f}+\frac{5\psi'''^2}{6\psi''^2}-\frac{\psi^{iv}}{2\psi''}\br\}\]
with all derivatives being evaluated at $x=k\pi$. This yields
\[c_2=\frac{1}{4}\bl\{2(1+a^2)(k\pi)^2-12a k\pi+9\br\}.\]

Summing over all the peaks we then obtain
\[K_n\sim\pi\sqrt{\frac{\pi}{n}}\bl\{\sigma_1+\frac{1}{8n}\bl(2\pi^2(1+a^2)\sigma_3-12\pi a\sigma_2+9\sigma_1\br)\br\},\]
where
\[\sigma_m:=\sum_{k=1}^\infty k^m e^{-k\pi a}.\]
We have
\[\sigma_1=\frac{1}{4\sinh^2 \fs\pi a},\quad \sigma_2=\frac{\cosh \fs\pi a}{4\sinh^3 \fs\pi a}, \quad \sigma_3=\frac{2+\cosh \pi a}{8\sinh^4 \fs\pi a}.\]
Hence we obtain our final estimate in the form
\bee\label{e32}
K_n\sim \frac{\pi^{3/2}}{4n^{1/2}} \bl\{1+\frac{T_1}{8n}\br\}\mbox{cosech}^2 \fs\pi a\qquad (n\to\infty),
\ee
where
\[T_1=9-12\pi a\,\mbox{coth} \fs\pi a+\pi^2(1+a^2)\frac{(2+\cosh \fs\pi a)}{\sinh^2 \fs\pi a}\]
with $a>0$ fixed and of O(1).

In Table 2 we show computed values of $K_n$ compared with the asymptotic estimate (\ref{e32}) for different values of $n$ and the parameter $a$. It is seen that the agreement is quite good and improves with increasing $n$. However, since $\psi''(k\pi)$ scales like $k^{-2}$, the peaks progressively broaden as $k$ increases with the consequence that
the saddle-point approximation eventually breaks down. In addition, the parameter $a$ cannot be too small on account of
the fact that the envelope of the minima of the integrand, given by $e^{-ax}(1-1/x^2)^n$, presents a maximum value at
$x\simeq (2n/a)^{1/3}$ equal to approximately $\exp\,[-\f{3}{2}(2na^2)^{1/3}]$. We require this last quantity to be small for the satisfactory estimation of each peak. This results in the condition $a\gg (2n)^{-1/2}$.
\begin{table}[h]
\caption{\footnotesize{Values of $K_n$ compared with asymptotic estimate (\ref{e32}).}}
\begin{center}
\begin{tabular}{|r|l|l||l|l|}
\hline
\mcol{1}{|c|}{} & \mcol{2}{c||}{$a=1$} & \mcol{2}{c|}{$a=3/2$}\\
\mcol{1}{|c|}{$n$} & \mcol{1}{c|}{$K_n$} & \mcol{1}{c||}{Asymptotic} & \mcol{1}{c|}{$K_n$} & \mcol{1}{c|}{Asymptotic} \\
[.1cm]\hline
&&&&\\[-0.25cm]
100 & 0.02707847 & 0.02689533 & 0.00523230 & 0.00521489 \\
200 & 0.01884203 & 0.01880232 & 0.00364706 & 0.00364449 \\
500 & 0.01181371 & 0.01180983 & 0.00228888 & 0.00228866 \\
1000& 0.00833214 & 0.00833153 & 0.00161452 & 0.00161448 \\
2000& 0.00588457 & 0.00588447 & 0.00114026 & 0.00114025 \\
4000& 0.00415855 & 0.00415854 & 0.00080580 & 0.00080580 \\
[.1cm]\hline
\mcol{1}{|c|}{} & \mcol{2}{c||}{$a=1/2$} & \mcol{2}{c|}{$a=2$}\\
\mcol{1}{|c|}{$n$} & \mcol{1}{c|}{$K_n$} & \mcol{1}{c||}{Asymptotic} & \mcol{1}{c|}{$K_n$} & \mcol{1}{c|}{Asymptotic} \\
[.1cm]\hline
&&&&\\[-0.25cm]
100 & 0.19606514 & 0.19692975 & 0.00108887 & 0.00108697 \\
200 & 0.13567443 & 0.13484945 & 0.00075359 & 0.00075332 \\
500 & 0.08386120 & 0.08361625 & 0.00047067 & 0.00047064 \\
1000& 0.05878333 & 0.05873199 & 0.00033143 & 0.00033143 \\
2000& 0.04139902 & 0.04139062 & 0.00023387 & 0.00023387 \\
4000& 0.02921970 & 0.02921838 & 0.00016520 & 0.00016520 \\
[.1cm]\hline
\end{tabular}
\end{center}
\end{table}

A closely related integral is
\[{\hat K}_n=\int_1^\infty e^{-ax}\bl(1-\frac{\cos^2 x}{x^2}\br)^{\!n}dx.\]
The peaks in the graph of the integrand are similar to those indicated in Fig.~1 but now occur at $x=(k+\fs)\pi$, $k=0, 1, 2, \ldots\ $. The lower limit of integration is chosen  and to lie in the interval $(\delta,\fs\pi-\delta)$ (with $\delta>0$ so as to avoid the origin and $\fs\pi$).
With $\psi(x)=-\log\,(1-\cos^2x/x^2)$ we find 
\[\psi''((k+\fs)\pi)=\frac{2}{(k+\fs)^2\pi^2}.\]
Then by similar arguments we obtain the leading asymptotic approximation
\[{\hat K}_n\sim \pi e^{-\pi a/2} \sqrt{\frac{\pi}{n}} \sum_{k=0}^\infty (k+\fs)e^{-k\pi a}
=\frac{\pi^{3/2} \cosh \fs\pi a}{4n^{1/2} \sinh^2 \fs\pi a}\qquad (n\to\infty).\]

\vspace{0.6cm}

\begin{center}
{\bf 4.\ An asymptotic expansion for Ball's integral}
\end{center}
\setcounter{section}{4}
\setcounter{equation}{0}
\renewcommand{\theequation}{\arabic{section}.\arabic{equation}}
Ball's integral is given by \cite{KB}
\bee\label{e41}
L(\nu;n)=\int_0^\infty \bl(\frac{\g(1+\nu) |J_\nu(x)|}{(x/2)^\nu}\br)^{\!n} x^{2\nu-1}dx,\qquad n\geq 2,\ \nu\geq \fs,
\ee
where $J_\nu(x)$ is the Bessel function of the first kind and $n$ is not necessarily an integer. In the case $\nu=\fs$
and integer values of $n$ the integral (\ref{e41}) (without the modulus signs) reduces to that in (\ref{e11}), since $\sqrt{\pi/(2x)} J_{1/2}(x)=\sin x/x$. The expansion of $L(\nu;n)$ for $n\to\infty$ has been derived in \cite{KOS}; here we revisit this result by means of the transformation used in Section 1.

We divide the integration path into the intervals $[0,j_{\nu,1}]$ and $[j_{\nu,1},\infty)$, where $j_{\nu,1}$ denotes the first zero of $J_\nu(x)$. In $[0,j_{\nu,1}]$, we have $J_\nu(x)\geq0$ and the modulus signs may be dropped in this interval, where
\[\sigma(x):=\frac{\g(1+\nu) J_\nu(x)}{(x/2)^\nu}=\sum_{k=0}^\infty \frac{(-x^2/4)^k}{k! (\nu+1)_k}\]
with $(a)_k=\g(a+k)/\g(a)$ being Pochhammer's symbol.
We set $\psi(x)=-\log\,\sigma(x)$ in the interval $[0,j_{\nu,1}]$, 
so that the integral (\ref{e41}) becomes
\bee\label{e42}
L(\nu;n)=\int_0^{j_{\nu,1}} e^{-n\psi(x)} x^{2\nu-1}dx+R_n(x),
\ee
where
\[R_n(x)=\int_{j_{\nu,1}}^\infty \bl(\frac{\g(1+\nu) |J_\nu(x)|}{(x/2)^\nu}\br)^{\!n} x^{2\nu-1}dx.
\]

The tail of the integral $R_n(x)$ satisfies the bound
\[|R_n(x)|<(2^\nu \g(1+\nu))^n \int_{j_{\nu,1}}^\infty \frac{dx}{x^{(n-2)\nu+1}}=\bl(\frac{2^\nu \g(1+\nu)}{j_{\nu,1}^\nu}\br)^n\,\frac{j_{\nu,1}^2}{(n-2)\nu}\]
since $|J_\nu(x)|\leq |J_\nu(j'_{\nu,2})|<1$ in $[j_{\nu,1},\infty)$, where $j'_{\nu,2}$ is the second zero of $J_\nu'(x)$. Defining the quantity 
\[\xi(\nu):=\frac{2^\nu \g(1+\nu)}{j_{\nu,1}^\nu},\]
and noting that $j_{\nu,1}=\pi$ when $\nu=\fs$, we see that $\xi(\fs)=2^{-1/2}$. Use of Stirling's approximation for the gamma function and the fact that \cite[p.~485]{WBF} $j_{\nu,1}>\nu$, shows that $\xi(\nu)\sim (2/e)^\nu \sqrt{2\pi\nu}\to 0$ as $\nu\to\infty$.
A plot of $\xi(\nu)$ for $\nu\geq\fs$ is shown in Fig. 2 where it is seen that $\xi(\nu)$ decreases monotonically\footnote{It is found that $\xi(\nu)$ is monotonically decreasing on $[0,\infty)$  with $\xi(\nu)<1$ for $\nu>0$.} with increasing $\nu$. Thus as $n\to\infty$, the bound on $R_n(x)$ when $\nu\geq\fs$ is of $O(n^{-1} 2^{-n/2})$ and so is exponentially small.
\begin{figure}[th]
	\begin{center} \includegraphics[width=0.5\textwidth]{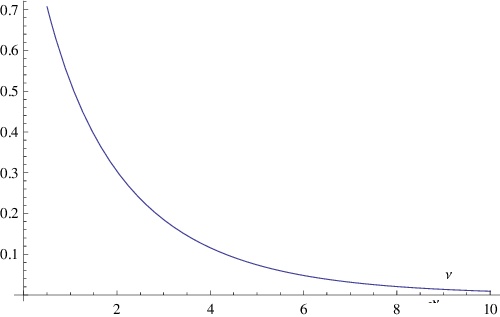}
\caption{\small{Plot of $\xi(\nu)$ against $\nu\geq\fs$.}} 
	\end{center}
\end{figure}

We now deal with the integral in (\ref{e42}), where we note that $\psi(0)=0$ and $\psi(j_{\nu,1})=\infty$.
Making the substitution $\tau^2=\psi(x)$, we obtain
\bee\label{e43}
L(\nu;n)=\int_0^\infty e^{-n\tau^2} x^{2\nu-1} \frac{dx}{d\tau}\,d\tau+R_n(x).
\ee
From the expansion
\[\tau^2=\psi(x)=\frac{x^2}{4(1+\nu)}+\frac{x^4}{32(1+\nu)^2(2+\nu)}+\frac{x^6}{96(1+\nu)^3(2+\nu)(3+\nu)}+\cdots\]
valid for $x<j_{\nu,1}$,
we find upon inversion
\[x=2(1+\nu)^{1/2} \bl\{\tau-\frac{\tau^3}{4(2+\nu)}-\frac{(1+11\nu)\tau^5}{96(2+\nu)^2(3+\nu)}-\frac{(17\nu^2-9\nu-20)\tau^7}{128(2+\nu)^3(3+\nu)(4+\nu)}+\cdots\br\}.\]
The last expansion holds in $\tau<\tau_0$, where $\tau_0^2=\log\,1/\sigma(j_{\nu+1,1})$, since $\psi'(x)=-J_{\nu+1}(x)/J_\nu(x)$ so that $x=j_{\nu+1,1}$ is the nearest point in the mapping $x\mapsto \tau$ where $dx/d\tau$ is singular (the point $x=j_{\nu,1}$ maps to $\infty$ in the $\tau$-plane).
This then yields the expansions
\[x^{2\nu-1}=(2(1+\nu)^{1/2}\tau)^{2\nu-1}\bl\{1-\frac{(2\nu-1)\tau^2}{4(2+\nu)}+\frac{(2\nu-1)(6\nu^2+\nu-19)\tau^4}{96(2+\nu)^2(3+\nu)}\hspace{3cm}\]
\[\hspace{6cm}-\frac{(2\nu-1)^2(2\nu^3-\nu^2-17\nu-20)\tau^6}{384(2+\nu)^3(3+\nu)(4+\nu)}+\cdots\br\}\]
and
\[\frac{dx}{d\tau}=2(1+\nu)^{1/2}\bl\{1-\frac{3\tau^2}{2(2+\nu)}-\frac{5(1+11\nu)\tau^4}{48(2+\nu)^2(3+\nu)}-\frac{7(17\nu^2-9\nu-20)\tau^6}{128(2+\nu)^3(3+\nu)(4+\nu)}+\cdots\br\}.\]

Combination of these last two expansions then produces
\bee\label{e44}
x^{2\nu-1} \frac{dx}{d\tau}=2^{2\nu}(1+\nu)^\nu \tau^{2\nu-1}\sum_{k=0}^\infty (-)^k b_k \tau^{2k}\qquad (\tau<\tau_0),
\ee
where 
\[b_0=1,\quad b_1=\frac{1+\nu}{2(2+\nu)},\quad b_2=\frac{3\nu^2+2\nu-5}{24(2+\nu)(3+\nu)},
\quad b_3=\frac{(1+\nu)(\nu^3-\nu^2-4\nu-8)}{48(2+\nu)^3(4+\nu)},\]
\[b_4=\frac{15\nu^7+15\nu^6-220\nu^5-918\nu^4+763\nu^3+15055\nu^2+26898\nu+13688}{5760(2+\nu)^3(3+\nu)^2(5+\nu)}, \ldots\ .\]
Insertion of the expansion (\ref{e44}) into the integral in (\ref{e43}) yields
\[2^{2\nu}(1+\nu)^\nu\int_0^\infty e^{-n\tau^2}\tau^{2\nu-1} \sum_{k=0}^\infty (-)^k b_k \tau^{2k} d\tau
\sim 2^{2\nu-1}(1+\nu)^\nu \g(\nu) \sum_{k=0}^\infty \frac{(-)^k b_k (\nu)_k}{n^{k+\nu}}.\]
As in Section 1, this will be a divergent expansion as we have integrated beyond the circle of convergence of $x^{2\nu-1} dx/d\tau$.

Thus, neglecting exponentially small terms we finally obtain the expansion
\bee\label{e45}
L(\nu;n)\sim 2^{2\nu-1}(1+\nu)^\nu \g(\nu) \sum_{k=0}^\infty \frac{(-)^k c_k}{n^{k+\nu}} \qquad (n\to\infty),
\ee
where the first few coefficients $c_k:=b_k(\nu)_k$ are
\[c_0=1,\quad c_1=\frac{\nu(1+\nu)\wp_1(\nu)}{2(2+\nu)},\quad c_2=\frac{\nu(1+\nu)\wp_2(\nu)}{24(2+\nu)(3+\nu)},\quad
c_3=\frac{\nu(1+\nu)\wp_3(\nu)}{48(2+\nu)^2(4+\nu)},\]
\[c_4=\frac{\nu (1 + \nu)\wp_4(\nu)}
{5760 (2 + \nu)^3 (3 + \nu) (5 + \nu)},\quad
c_5=\frac{\nu (1 + \nu)\wp_5(\nu) }{
11520 (2 + \nu)^4 (3 + \nu) (6 + \nu)},\] 
\[c_6=      
      \frac{\nu (1 + 
      \nu)\wp_6(\nu)}{2903040 (2 + \nu)^5 (3 + \nu)^2 (4 + \nu) (7 + 
      \nu)} ,\]
with the polynomials $\wp_k(\nu)$  given by
\begin{eqnarray*}     
\wp_1(\nu)\!\!&=&\!\!1,\quad \wp_2(\nu)=3\nu^2+2\nu-5,\quad
\wp_3(\nu)=(1+\nu)(\nu^3-\nu^2-4\nu-8),\\
\wp_4(\nu)\!\!&=&\!\!15 \nu^7+ 15 \nu^6 -220 \nu^5-918 \nu^4+ 763 \nu^3 + 15055 \nu^2 + 26898 \nu + 13688,\\
\wp_5(\nu)\!\!&=&\!\!3 \nu^9\!-\! 7 \nu^8\!-\! 66 \nu^7\!- \!246 \nu^6\! +\! 2307 \nu^5\! + 
     \! 6825 \nu^4\! -\! 43668 \nu^3\! -118508 \nu^2\!  -\! 89904 \nu\!  -\!19392,\\
\wp_6(\nu)\!\!&=&\!\!63 \nu^{13}\!-\! 3276 \nu^{11}\!-\! 16856 \nu^{10}\! +\! 131726 \nu^9\!+\! 781856 \nu^8\!-\!
      4685840 \nu^7 \!-\! 14835768 \nu^6\\
 &&  \!\!\!\!    + 104879595 \nu^5\! \!+ 322760624 \nu^4\!-\! 
      328990364 \nu^3\!  -\! 1748824256 \nu^2\!-\! 1801386304 \nu\!-\!590749440 .
      \end{eqnarray*}
The first coefficients (with $k\leq 3$) agree with those found by Kerman {\it et al.} \cite{KOS}. It is easily verified that when $\nu=\fs$ the above coefficients agree with those listed in Table 1.

The estimate for the tail of the integral 
over the interval $[2^\nu \g(1+\nu),\infty)$ in \cite{KOS}, however, is only $O(n^{-1})$, which is not sufficiently sharp to justify the expansion of the main integral beyond its leading term. We have demonstrated that the tail of the integral (\ref{e41}) is {\it exponentially small} as $n\to\infty$. 

An obvious extension of (\ref{e41}) is the integral
\[{\cal L}(\nu,a;n)=\int_0^\infty \bl(\frac{\g(1+\nu) |J_\nu(x)|}{(x/2)^\nu}\br)^{\!n} x^{a-1}dx, \qquad a>0,\]
where $a$ is fixed and it is assumed that $n$ satisfies the condition $n(\nu+\fs)>a$ to secure convergence at infinity.
The procedure described above then produces the expansion (when exponentially small terms are neglected)
\bee\label{e46}
{\cal L}(\nu,a;n)\sim 2^{a-1} (1+\nu)^{a/2} \g(\fs a) \sum_{k=0}^\infty \frac{(-)^k d_k}{n^{k+a/2}}\qquad (n\to\infty),
\ee
where the first few coefficients $d_k$ are given by
\[d_0=1, \quad d_1=\frac{(\fs a)_2}{2(2+\nu)},\quad d_2=\frac{(\fs a)_3}{48(2+\nu)^2(3+\nu)}((3a-14)\nu+9a-10),\]
\[d_3=\frac{(\fs a)_4}{192(2+\nu)^3(3+\nu)(4+\nu)}\bl\{(a^2-14a+64)\nu^2+(7a^2-66a+32)\nu+4(a-4)(3a+2)\br\},\]
\[d_4=\frac{(\fs a)_5}{46080(2+\nu)^4(3+\nu)^2(4+\nu)(5+\nu)}\bl\{(15a^3-420a^2+4820a- 23824) \nu^4\hspace{2.4cm}\] \[+(225a^3-5340a^2+42860a- 65776)\nu^3 + (1245a^3-23340a^2+103740a+100560)\nu^2\]
\[ +(3015a^3-39300a^2+45940a + 252784)\nu 
 + 2700 a^3- 18000 a^2 - 18800 a+109504\br\}.\]  
When $a=2\nu$, it is seen that the $d_k$ reduce to the coefficients $c_k$ ($k\leq 4$) appearing in the expansion  (\ref{e45}).

\begin{table}[th]
\caption{\footnotesize{Values of the absolute relative error in the computation of ${\cal L}(\nu,a;n)$ against truncation index $k$ when $\nu=4/3$ and $n=100$.}}
\begin{center}
\begin{tabular}{|r|l|l|l|}
\hline
\mcol{1}{|c|}{$k$} & \mcol{1}{c|}{$a=8/3$} & \mcol{1}{c|}{$a=2/3$} & \mcol{1}{c|}{$a=10/3$}\\
[.1cm]\hline
&&&\\[-0.25cm]
0 & $4.664\times 10^{-03}$ & $6.676\times 10^{-04}$ & $6.565\times10^{-03}$ \\
1 & $2.738\times 10^{-06}$ & $8.987\times 10^{-07}$ & $1.047\times10^{-05}$ \\
2 & $3.307\times 10^{-08}$ & $6.661\times 10^{-10}$ & $6.041\times10^{-08}$ \\
3 & $4.006\times 10^{-10}$ & $2.405\times 10^{-11}$ & $5.961\times10^{-10}$ \\
4 & $2.914\times 10^{-12}$ & $3.655\times 10^{-13}$ & $2.743\times10^{-12}$ \\
[.1cm]\hline
\end{tabular}
\end{center}
\end{table}

In Table 2 we present the values of the absolute relative error in the evaluation of ${\cal L}(\nu,a;n)$ using the expansion (\ref{e46}) for different truncation index $k$. The first column shows the values $\nu=4/3$, $a=8/3$, which corresponds to the integral (\ref{e41}).

\end{document}